\documentclass[12pt]{article}
\usepackage{amsmath,amsthm,amssymb}
\usepackage{mathrsfs}
\usepackage[titletoc]{appendix}
\usepackage{bbm}
\usepackage{float}
\usepackage{tikz}
\begin{document}
\input amssym.def
\newcommand{\singlespace}{
    \renewcommand{\baselinestretch}{1}
\large\normalsize}
\newcommand{\doublespace}{
   \renewcommand{\baselinestretch}{1.2}
   \large\normalsize}
\renewcommand{\theequation}{\thesection.\arabic{equation}}

\setcounter{equation}{0}
\def \ten#1{_{{}_{\scriptstyle#1}}}
\def \Z{\Bbb Z}
\def \C{\Bbb C}
\def \R{\Bbb R}
\def \Q{\Bbb Q}
\def \N{\Bbb N}
\def \l{\lambda}
\def \V{V^{\natural}}
\def \wt{{\rm wt}}
\def \tr{{\rm tr}}
\def \Res{{\rm Res}}
\def \End{{\rm End}}
\def \Aut{{\rm Aut}}
\def \mod{{\rm mod}}
\def \Hom{{\rm Hom}}
\def \im{{\rm im}}
\def \<{\langle}
\def \>{\rangle}
\def \w{\omega}
\def \c{{\tilde{c}}}
\def \o{\omega}
\def \t{\tau }
\def \ch{{\rm ch}}
\def \a{\alpha }
\def \b{\beta}
\def \e{\epsilon }
\def \la{\lambda }
\def \om{\omega }
\def \O{\Omega}
\def \qed{\mbox{ $\square$}}
\def \pf{\noindent {\bf Proof: \,}}
\def \voa{vertex operator algebra\ }
\def \voas{vertex operator algebras\ }
\def \p{\partial}
\def \1{{\bf 1}}
\def \ll{{\tilde{\lambda}}}
\def \H{{\bf H}}
\def \F{{\bf F}}
\def \h{{\frak h}}
\def \g{{\frak g}}
\def \rank{{\rm rank}}
\def \({{\rm (}}
\def \){{\rm )}}
\def \Y {\mathcal{Y}}
\def \I {\mathcal{I}}
\def \A {\mathcal{A}}
\def \B {\mathcal {B}}
\def \Cc {\mathcal {C}}
\def \H {\mathcal{H}}
\def \M {\mathcal{M}}
\def \V {\mathcal{V}}
\def \O{{\bf O}}
\def \AA{{\bf A}}
\singlespace
\newtheorem{thm}{Theorem}[section]
\newtheorem{prop}[thm]{Proposition}
\newtheorem{lem}[thm]{Lemma}
\newtheorem{cor}[thm]{Corollary}
\newtheorem{rem}[thm]{Remark}
\newtheorem*{CPM}{Theorem}
\newtheorem{definition}[thm]{Definition}

\begin{center}
{\Large {\bf  Bimodule and twisted representation of vertex operator algebras }} \\

\vspace{0.5cm} Qifen Jiang\footnote{Supported by NSFC grants 11101269 and 11431010.\\
 \indent E-mail:qfjiang@sjtu.edu.cn}
\\
Department of Mathematics, Shanghai Jiaotong University\\ Shanghai 200240, China  \\
Xiangyu Jiao\footnote{Supported by NSFC grant 11401213\\
\ \indent E-mail:xyjiao@math.ecnu.edu.cn}\\
Department of Mathematics, East China Normal University\\ Shanghai 200241, China \\
\end{center}
\hspace{1.5 cm}

{\bf Abstract}
In this paper, for a vertex operator algebra $V$ with an automorphism $g$ of order $T,$ an admissible $V$-module $M$ and a fixed nonnegative rational number $n\in\frac{1}{T}\Z_+,$ we construct an $A_{g,n}(V)$-bimodule
$\AA_{g,n}(M)$ and study its some properties, discuss the connections between bimodule $\AA_{g,n}(M)$ and intertwining operators. Especially, bimodule $\AA_{g,n-\frac{1}{T}}(M)$ is a natural quotient of $\AA_{g,n}(M)$ and there is a linear isomorphism between the space ${\cal I}_{M\,M^j}^{M^k}$ of intertwining operators and the space of homomorphisms  $\Hom_{A_{g,n}(V)}(\AA_{g,n}(M)\otimes_{A_{g,n}(V)}M^j(s), M^k(t))$ for $s,t\leq n, M^j, M^k$ are $g$-twisted $V$ modules, if $V$ is $g$-rational.


{\bf Remark:} It has been accepted for publication in SCIENCE CHINA Mathematics.
\section{Introduction}

 \hspace{0.5cm} It is known that twisted representations are the main ingredients in orbifold conformal field theory (refer [DHVW1-DHVW2, LJ1-LJ2, FLM, DM, DLM1, DLM3, DLM5, DLM6 and MT]etc). They play a fundamental role in the construction of the moonshine vertex operator algebra $V^{\natural }$ \cite{FLM} and other orbifold vertex operator algebras \cite{DGM}. Many of them have been studied in the literature (refer [ABD, DJ2, DLM3, DJM4, HY, MT]). In \cite{DLM3, Z}, associative algebra $A_{g}(V)$ was constructed for a given vertex operator algebra $V$ with an automorphism $g$ of finite order $T$. This is helpful to investigate abstract orbifold models. $A_{g,n}(V),$ which was constructed in \cite{DLM5} for nonnegative numbers $n\in\frac{1}{T}\Z_{+},$  is a generalization of the algebra $A_g(V),$ where $A_g(V)=A_{g,0}(V)$ (where $A_{id,0}(V)$ coincides with $A(V)$ in \cite{Z}, and $A_{id,n}(V)$ is exactly $A_n(V)$ in \cite{DLM4}). It was shown that $A_{g,n-\frac{1}{T}}(V)$ is a natural quotient of $A_{g,n}(V)$.

On the other hand, an $A_{n}(V)$-bimodule $\AA_{n}(M)$ for any $V$-module $M$ was also introduced in \cite{FZ, DR} to deal with the intertwining operators and fusion rules.
Moreover an important result in \cite{FZ, L2, DR} is that
there is an isomorphism between the space of intertwining operators among irreducible $V$-modules and a certain space constructed from bimodules for  $A_{id,n}(V)$. By generalizing the bimodule theory in \cite{FZ}, a sequence of bimodules of $A_{g,n}(V)$-$A_{g,m}(V)$ in \cite{DJ1, DJ2} were studied, which were denoted by $A_{g,n,m}(V)$.

Motivated by \cite{DR, DJ1, DJ2},
for any $g,\ n$ as above, in this paper we construct an $A_{g,n}(V)$-bimodule
$\AA_{g,n}(M)$ for any admissible $V$-module $M$ such that $\AA_{g,n-\frac{1}{T}}(M)$ is a natural quotient of $\AA_{g,n}(M).$ The structure of those bimodules are discussed. We also consider the connections between bimodule $\AA_{g,n}(M)$ and intertwining operators. Moreover, it is established that if $V$ is $g$-rational then there is a linear isomorphism between the space ${\cal I}_{M\,M^j}^{M^k}$ of intertwining operators and the space of homomorphism $$\Hom_{A_{g,n}(V)}(\AA_{g,n}(M)\otimes_{A_{g,n}(V)}M^j(s), M^k(t))$$ for $s,t\leq n, M^j, M^k$ are $g$-twisted $V$ modules.

The paper is organized as follows. We recall the construction of associative algebras $A_{g,n}(V)$ and some results
from \cite{DLM3}-\cite{DLM5} in Section 2. In Section 3,we give the construction of $\AA_{g,n}(M)$ and
study its  structure.
In section 4, we discuss the relation between $\AA_{g,n}(M)$ and intertwining operators. As in \cite{FZ, L2, DR},  we  obtain the isomorphism from the space ${\cal I}_{M\,M^j}^{M^k}$ of intertwining operators to $\Hom_{A_{g,n}(V)}(\AA_{g,n}(M)\otimes_{A_{g,n}(V)}M^j(s), M^k(t))$
if $V$ is $g$-rational.

\section{The associative algebras $A_{g,n}(V)$}
\setcounter{equation}{0}

\hspace{0.5cm}
Let $V=(V,Y,{\bold 1},\omega)$ be a vertex operator algebra (\cite{B}, \cite{FLM}, \cite{LL}) with an automorphism $g$ of finite order $T$.
Then $V$ has the decomposition of eigenspace with respect to the action $g$:
$$V =\bigoplus_{r\in\Z/{T\Z}} V^r$$where $V^r=\{v\in V|gv=e^{-\frac{2\pi ir}{T}}v\}.$
 Firstly, we recall definitions of modules of different types from \cite{DLM3, FLM, FFR, D}.
\begin{definition} A {\em weak g-twisted V-module} $M$ is a vector space
equipped with a linear map
\begin{equation*}
\begin{split}
Y_M(\cdot, z) : &V \to (\End M)[[z^{\frac{1}{T}}, z^{-\frac{1}{T}}]]\\
&v  \mapsto Y_M(v, z) =\sum_{n\in\Q}v_nz^{-n-1}~~~~~~(v_n \in \hbox{End}M)
\end{split}
\end{equation*}
which satisfies the following conditions for all $0\leq r\leq T-1, u \in V^r,$ $v \in V,$ $w \in M$.
\begin{equation*}
\begin{split}
&Y_M(u, z)=\sum_{n\in\frac{r}{T}+\Z}u_nz^{-n-1};\\
&u_nw = 0 \hbox{ for } n \gg 0;\\
& Y_M(\1, z) = \hbox{id}_M;\\
z_0^{-1}\delta(&\frac{z_1-z_2}{z_0})Y_M(u, z_1)Y_M(v, z_2)-z_0^{-1}\delta(\frac{z_2-z_1}{-z_0})Y_M(v, z_2)Y_M(u, z_1)\\
&= z_2^{-1}\delta(\frac{z_1-z_0}{z_2})(\frac{z_1-z_0}{z_2})^{-\frac{r}{T}}Y_M(Y (u, z_0)v, z_2).
\end{split}
\end{equation*}
\end{definition}

\begin{definition} An ({\em ordinary}) g-twisted $V$-module is a weak g-twisted $V$-module $M$ which carries a $\C$-grading induced by the spectrum
of $L(0)$ where $L(0)$ is a component operator of
$$Y_M(\omega, z) =\sum_{n\in\Z}L(n)z^{-n-2}.$$
That is$$M =\oplus _{\l\in \C}M_{\l},$$ where $M_\l = \{w\in M|L(0)w = \l w\}$. Moreover one requires that $M_\l$ is
finite dimensional and for fixed $\l$, $M_{\frac{n}{T}+\l}= 0$ for all small enough integers $n$.
\end{definition}

\begin{definition} An {\em admissible} g-twisted $V$-module is a weak g-twisted $V$-module $M$ which carries a $\frac{1}{T}\Z_+$-grading
$$M =\oplus_{n\in\frac{1}{T}\Z_+} M(n)$$ that satisfies the following
$$v_mM(n) \subset M(n+\wt v-m-1)$$
for homogeneous $v\in V.$
\end{definition}

\begin{rem}\label{rem1} Assume $M$ is an irreducible addimisible $g$-twisted $V$-module, one knows that $M=\oplus_{n\in \frac{1}{T}\Z_+} M_{n+h}$ for some $h\in \C$ such that $M_h\neq 0,$  where $h$ is the conformal weight of $M.$ We always assume $M(n)=M_{n+h}.$ \end{rem}
Next we recall some results of the associative algebra $A_{g,n}(V)$ from [DLM3-DLM5]. Also see \cite{Z}.

Fix $n=l+\frac{i}{T}\in\frac{1}{T}\Z_+$ with $l$ a nonnegative integer and $0\leq i\leq T-1$. For $0\leq r\leq T-1$, define
 $$
\delta_i(r)=\begin{cases} &1,\quad\text{if}\ i\geq r,\\
&0,\quad\text{if}\ i< r.
\end{cases}
.$$Let $O_{g,n}(V)$ be the linear span of all $u\circ_{g,n} v$ and $L(-1)u+L(0)u$
where for homogeneous $u\in V^r$ and $v\in V,$
\begin{equation*}\label{2.1}
u\circ_{g,n} v=\Res_{z}Y(u,z)v\frac{(1+z)^{\wt u+l-1+\delta_i(r)+\frac{r}{T}}}{z^{2l+\delta_i(r)+\delta_i(T-r)+1}}.
\end{equation*}
We also define a second product $*_{g,n}$ on $V$ for $u\in V^r$ and $v\in V$ as
follows:
\begin{eqnarray*}\label{2.2}
& & u*_{g,n}v=\sum_{m=0}^{l}(-1)^m{m+l\choose l}\Res_zY(u,z)\frac{(1+z)^{\wt\,u+l}}{z^{l+m+1}}v\\
& &\ \ \ \ \ \ \ \ =\sum_{m=0}^l\sum_{j=0}^{\infty}(-1)^m
{m+l\choose l}{\wt u+l\choose j}u_{j-m-l-1}v.\nonumber
\end{eqnarray*}
if $r=0$ and $u*_{g,n}v=0$ if $r>0$.

Define $A_{g,n}(V)$ to be the quotient $V/O_{g,n}(V).$ Then $A_{g,0}(V)=A_{g}(V)$ and $A_{id,n}(V)=A_{n}(V)$ have already been defined and studied in \cite{DLM3} and \cite{DLM4} respectively.

The following theorem summarizes the main results of \cite{DLM5}.
\begin{thm}\label{tha} Let $V$ be a vertex operator algebra and $g$ an automorphism of $V$ of finite order $T$. Let
$M =\oplus_{m\in\frac{1}{T}\Z_+} M(m)$ be an admissible $g$-twisted $V$-module. Let $n\in\frac{1}{T}\Z_+$. Then

 $(1)$ $A_{g,n}(V)$ is an associative algebra whose product is induced by
$*_{g,n}.$

$(2)$ The identity map on $V$ induces an algebra epimorphism from
$A_{g,n}(V)$ to $A_{g,n-1/T}(V).$

$(3)$ Let $W$ be a weak $g$-twisted $V$ module and set
$$\Omega_n(W)=\{w\in W|u_{wt u-1+k}w=0, u\in V, k>n\}.$$
Then $\Omega_n(W)$ is an $A_{g,n}(V)$-module such that $v+O_{g,n}(V)$ acts as $o(v)=v_{\wt v-1}$
 for homogeneous $v$.

$(4)$ Each $M(m)$ for $m\leq n$ is an $A_{g,n}(V)$-submodule of $\Omega_n(W).$ Furthermore,
$M$ is irreducible if and only if each $M(m)$ is an irreducible
$A_{g,n}(V)$-module.

$(5)$ The linear map $\varphi$ from $V$ to $V$ defined by $\varphi(u)=e^{L(1)}(-1)^{L(0)}u$ for $u\in V$ induces an anti-isomorphism from $A_{g,n}(V)$ to $A_{g^{-1},n}(V)$.
\end{thm}

\section{Bimodules $\AA_{g,n}(M)$ of $A_{g,n}(V)$}

\hspace{0.5cm}Motivated by the ideas of $A(V)$-bimodule $\AA(M)$ and $A_n(V)$-bimodule $\AA_n(M)$  from \cite{Z} and \cite{DR} for any nonnegative integer $n$ respectively, we will construct and study $A_{g,n}(V)$-bimodule $\AA_{g,n}(M)$ for $n\in \frac{1}{T}\Z_+$ in this section.

Let $V$ be same as vertex operator algebra in section 2 and $M$ be an admissible $V$-module, $\O_{g,n}(M)$ be the linear span of $u\circ_{g,n} w$
where for homogeneous $u\in V^r$ and $w\in M,$
\begin{equation*}\label{2.1m}
u\circ_{g,n} w=\Res_{z}Y_M(u,z)w\frac{(1+z)^{\wt u+l-1+\delta_i(r)+\frac{r}{T}}}{z^{2l+\delta_i(r)+\delta_i(T-r)+1}}.
\end{equation*}
We also define a left bilinear product $*_{g,n}$  for $u\in V^r$ and $w\in M:$
\begin{eqnarray}\label{2.2m}
& & u*_{g,n}w=\sum_{m=0}^{l}(-1)^m{m+l\choose l}\Res_zY_M(u,z)\frac{(1+z)^{\wt\,u+l}}{z^{l+m+1}}w\\
& &\ \ \ \ \ \ \ \ =\sum_{m=0}^l\sum_{j=0}^{\infty}(-1)^m
{m+l\choose l}{\wt u+l\choose j}u_{j-m-l-1}w,\nonumber
\end{eqnarray}
if $r=0$ and $u*_{g,n}w=0$ if $r>0$
and a right bilinear  product
\begin{eqnarray}\label{2.3m}
& & w*_{g,n}u=\sum_{m=0}^{l}(-1)^l{m+l\choose l}\Res_zY_M(u,z)\frac{(1+z)^{\wt\,u+m-1}}{z^{l+m+1}}w\\
& &\ \ \ \ \ \ \ \ =\sum_{m=0}^l\sum_{j=0}^{\infty}(-1)^l
{m+l\choose l}{\wt u+m-1\choose j}u_{j-m-l-1}w.\nonumber
\end{eqnarray}
if $r=0$ and $w*_{g,n}u=0$ if $r>0$.

Set $\AA_{g,n}(M)=M/\O_{g,n}(M).$ One can see that $\AA_{1,n}(M)$ is the $A_n(V)$-bimodule $\AA_n(M)$ studied in \cite{DR}. Moreover, $\AA_{1,0}(M)$ is the $A(V)$-bimodule $\AA(M)$ studied in \cite{FZ}.

\begin{lem}\label{l2.2} $(1)$
 Assume that $u\in V^r$ is homogeneous,
$w\in M$ and $m\ge k\ge 0.$ Then
$$\Res_{z}Y_M(u,z)w\frac{(1+z)^{{\wt}u+l-1+\delta_i(r)+\frac{r}{T}+k}}{z^{2l+\delta_i(r)+\delta_i(T-r)+1+m}}\in \O_{g,n}(M).$$

$(2)$ For homogeneous $u\in V^0$ and $w\in M,$
$$u*_{g,n}w-w*_{g,n}u-\Res_zY_M(u,z)w(1+z)^{\wt u-1}\in \O_{g,n}(M).$$
\end{lem}

\pf The proof of $(1)$ is similar to that of Lemma 2.1.2 of [Z].\\ 
For $(2)$, it follows directly
from the definitions (\ref{2.2m}) and (\ref{2.3m}) and using the following
result in the appendix of \cite{DLM4}: 
$$\sum_{m=0}^l{m+l\choose l}\frac{(-1)^m(1+z)^{l+1}-(-1)^l(1+z)^m}{z^{l+m+1}}=1.$$
Thus we can obtain 
the result in $(2)$.\qed

\begin{lem}\label{l2.f1}

 $$(L(-1)+L(0))V*_{g,n}M\subset \O_{g,n}(M), M*_{g,n} (L(-1)+L(0))V\subset \O_{g,n}(M)$$
\end{lem}
\begin{pf} By the definition of $*_{g,n}$, we only need to prove the case $u\in V^0, w\in M$. From (\ref{2.2m}), we have
\begin{eqnarray*}
&&(L(-1)u)*_{g,n}w=\sum_{m=0}^l(-1)^m{m+l\choose l}\Res_zY_M(L(-1)u,z)w\frac{(1+z)^{\wt\,u+1+l}}{z^{l+m+1}}\\
&&=\sum_{m=0}^l(-1)^{m+1}{m+l\choose l}\Res_zY_M(u,z)w\frac{d}{dz}(\frac{(1+z)^{\wt\,u+1+l}}{z^{l+m+1}})\end{eqnarray*}
\begin{eqnarray*}
&&=\sum_{m=0}^l(-1)^{m+1}{m+l\choose l}\Res_zY_M(u,z)w[\frac{(\wt\,u+1+l)(1+z)^{\wt\,u+l}}{z^{l+m+1}}\\
&&\ \ \ \ \ \ \ \ \ \ \ \ \ \ \ \ \ \ \ \ \ \ \ \ \ \ \ \ \ \ \ \ \ -\frac{(l+m+1)(1+z)^{\wt\,u+1+l}}{z^{l+m+2}}],\\
&&(L(0)u)*_{g,n}w=\sum_{m=0}^l(-1)^m{m+l\choose l}\Res_zY_M(L(0)u,z)w\frac{(1+z)^{\wt\,u+l}}{z^{l+m+1}}\\
&&=\sum_{m=0}^l(-1)^m{m+l\choose l}\wt\,u\Res_zY_M(u,z)w\frac{(1+z)^{\wt\,u+l}}{z^{l+m+1}}\end{eqnarray*}
Thus \begin{eqnarray*}&&(L(-1)+L(0))u*_{g,n}w\\
&&=\sum_{m=0}^l(-1)^{m}{m+l\choose l}\Res_zY_M(u,z)w[\frac{(-l-1)(1+z)^{\wt\,u+l}}{z^{l+m+1}}\\
&&\ \ \ \ \ \ \ \ \ \ \ \ \ \ \ \ \ \ \ \ \ \ \ \ \ \ \ \ \ \ \ \ \ +\frac{(l+m+1)(1+z)^{\wt\,u+1+l}}{z^{l+m+2}}]\\
&&=\sum_{m=0}^l(-1)^{m}{m+l\choose l}\Res_zY_M(u,z)w\frac{(mz+l+m+1)(1+z)^{\wt\,u+l}}{z^{l+m+2}}.\end{eqnarray*}
Using Lemma 2.2 in \cite{DLM4}:
$$\sum_{m=0}^l(-1)^{m}{m+l\choose l}\frac{(mz+l+m+1)(1+z)^{\wt\,u+l}}{z^{l+m+2}}=(-1)^{l}{2l+1\choose l}\frac{2l+1}{z^{2l+2}},$$
we obtain $$(L(-1)+L(0))u*_{g,n}w=(-1)^{l}{2l+1\choose l}\Res_zY_M(u,z)w\frac{(2l+1)(1+z)^{\wt\,u+l}}{z^{2l+2}}.$$
Since $u\in V^0, \delta_i(r)=1$ and $\delta_i(T-r)=0$, it is clear that the right-hand side lies in $\O_{g,n}(M)$.
So$$(L(-1)+L(0))V*_{g,n}M\subset \O_{g,n}(M)$$

For the second containment,we only need to use this result and Lemma \ref{l2.2}(2). This completes the proof.
\end{pf}\qed

\begin{lem}\label{l2.3}
$(1)\ O_{g,n}(V)*_{g,n}M\subset \O_{g,n}(M), M*_{g,n}O_{g,n}(V)\subset \O_{g,n}(M),$
$$(2)\ V*_{g,n}\O_{g,n}(M)\subset \O_{g,n}(M), \O_{g,n}(M)*_{g,n}V\subset \O_{g,n}(M).$$
\end{lem}

\begin{pf}  By Lemma \ref{l2.f1}, it remains to prove that
$(u\circ_{g,n}v)*_{g,n}w, w*_{g,n}(u\circ_{g,n} v), u*_{g,n}(v\circ_{g,n}w), (v\circ_{g,n}w)*_{g,n}u\in \O_{g,n}(M)$ for $u,v\in V$ and $w\in M.$

(1) Since $V^r*_{g,n}M$ and $M*_{g,n}V^r$ equal zero for $r\neq 0$, we only need to consider the case $u\in V^r$ and $v\in V^{T-r}$.
It is not difficult to prove $(u\circ_{g,n}v)*_{g,n}w\subset \O_{g,n}(M)$ according to Lemma 3.4(1) in \cite{DR}, we omit it. We now prove $w*_{g,n}(u\circ_{g,n} v)\in \O_{g,n}(M)$.

Let $u\in V^r, v\in V^{T-r}$ and $w\in M$, we have
$$(u\circ_{g,n}v)=\sum_{j\geq 0}{\wt u+l-1+\delta_i(r)+\frac{r}{T}\choose j}(u_{j-2l-1-\delta_i(r)-\delta_i(T-r)}v).$$
Then it follows Lemma \ref{l2.2}(2) that
\begin{eqnarray*}
&&w*_{g,n}(u\circ_{g,n} v)\\
&&\equiv -\sum_{j\geq 0}{\wt u+l-1+\delta_i(r)+\frac{r}{T}\choose j}\Res_{z_2}Y_M(u_{j-2l-1-\delta_i(r)-\delta_i(T-r)}v,z_2)w\times\\
&&\ \ \ \ \ \ \ \ \ \ \ \ \ \ (1+z_2)^{\wt u+\wt v+2l+\delta_i(r)+\delta_i(T-r)-1-j}(\mod \O_{g,n}(M))\\ \end{eqnarray*}
\begin{eqnarray*}
&&=-\Res_{z_2}\Res_{z_0}Y_M(Y(u,z_0)v,z_2)w\frac{(1+z_2+z_0)^{\wt u+l-1+\delta_i(r)+\frac{r}{T}}}{z_0^{2l+\delta_i(r)+\delta_i(T-r)+1}}\times\\
&&\quad\qquad\qquad (1+z_2)^{\wt v+l+\delta_i(T-r)-\frac{r}{T}}(\mod \O_{g,n}(M))\\
&&=-\Res_{z_1}\Res_{z_2}Y_M(u,z_1)Y_M(v,z_2)w\frac{(1+z_1)^{\wt u+l-1+\delta_i(r)+\frac{r}{T}}(1+z_2)^{\wt v+l+\delta_i(T-r)-\frac{r}{T}}}{(z_1-z_2)^{2l+\delta_i(r)+\delta_i(T-r)+1}}\\
&&+\Res_{z_2}\Res_{z_1}Y_M(v,z_2)Y_M(u,z_1)w
\frac{(1+z_1)^{\wt u+l-1+\delta_i(r)+\frac{r}{T}}(1+z_2)^{\wt v+l+\delta_i(T-r)-\frac{r}{T}}}{(-z_2+ z_1)^{2l+\delta_i(r)+\delta_i(T-r)+1}}\\
&&=\sum_{j\geq 0}(-1)^{j+1}{-2l-1-\delta_i(r)-\delta_i(T-r)\choose j}\Res_{z_1}\Res_{z_2}Y_M(u,z_1)Y_M(v,z_2)w\times\\
&&\qquad\qquad\frac{(1+z_1)^{\wt u+l-1+\delta_i(r)+\frac{r}{T}}(1+z_2)^{\wt v+l+\delta_i(T-r)-\frac{r}{T}}}{z_1^{2l+\delta_i(r)+\delta_i(T-r)+1+j}z_2^{-j}}\\
&&+\sum_{j\geq 0}(-1)^{\delta_i(r)+\delta_i(T-r)+1}{-2l-1-\delta_i(r)-\delta_i(T-r)\choose j}\Res_{z_2}\Res_{z_1}Y_M(v,z_2)Y_M(u,z_1)w\times\\
&&\qquad\qquad\frac{(1+z_1)^{\wt u+l-1+\delta_i(r)+\frac{r}{T}}(1+z_2)^{\wt v+l+\delta_i(T-r)-\frac{r}{T}}}{z_1^{-j}z_2^{2l+\delta_i(r)+\delta_i(T-r)+1+j}}.
\end{eqnarray*}
Since $u\in V^r, v\in V^{T-r}$ and from Lemma \ref{l2.2} we know that both
$$\Res_{z_1}\Res_{z_2}Y_M(u,z_1)Y_M(v,z_2)w\frac{(1+z_1)^{\wt u+l-1+\delta_i(r)+\frac{r}{T}}(1+z_2)^{\wt v+l+\delta_i(T-r)-\frac{r}{T}}}{z_1^{2l+\delta_i(r)+\delta_i(T-r)+1+j}z_2^{-j}}$$
and
$$\Res_{z_2}\Res_{z_1}Y_M(v,z_2)Y_M(u,z_1)w\frac{(1+z_1)^{\wt u+l-1+\delta_i(r)+\frac{r}{T}}(1+z_2)^{\wt v+l+\delta_i(T-r)-\frac{r}{T}}}{z_1^{-j}z_2^{2l+\delta_i(r)+\delta_i(T-r)+1+j}}$$
lie in $\O_{g,n}(M).$

$(2)$ Let $u\in V^0, v\in V^r$ and $w\in M$, we have
\begin{eqnarray*}
&&u*_{g,n}(v\circ_{g,n}w )
=\sum_{m=0}^l(-1)^m{m+l\choose l}\Res_{z_1}Y_M(u,z_1)\Res_{z_2}Y_M(v,z_2)w\\
&&\ \ \ \ \ \ \qquad\qquad\qquad\times\frac{(1+z_1)^{\wt u+l}}{z_1^{l+m+1}}\cdot\frac{(1+z_2)^{\wt v+l-1+\delta_i(r)+\frac{r}{T}}}{z_2^{2l+\delta_i(r)+\delta_i(T-r)+1}}\\
&&=\sum_{m=0}^l(-1)^m{m+l\choose l}\Res_{z_2}\Res_{z_1}Y_M(v,z_2)Y_M(u,z_1)w\\
&&\ \ \ \ \ \ \qquad\qquad\qquad\times\frac{(1+z_1)^{\wt u+l}}{z_1^{l+m+1}}\cdot\frac{(1+z_2)^{\wt v+l-1+\delta_i(r)+\frac{r}{T}}}{z_2^{2l+\delta_i(r)+\delta_i(T-r)+1}}\\
&&+\sum_{m=0}^l(-1)^m{m+l\choose l}\Res_{z_2}\Res_{z_0}Y_M(Y(u,z_0)v,z_2)w
\times\frac{(1+z_2+z_0)^{\wt u+l}}{(z_2+z_0)^{l+m+1}}\cdot\frac{(1+z_2)^{\wt v+l-1+\delta_i(r)+\frac{r}{T}}}{z_2^{2l+\delta_i(r)+\delta_i(T-r)+1}}\\\end{eqnarray*}
\begin{eqnarray*}
&&=\sum_{m=0}^l(-1)^m{m+l\choose l}\Res_{z_2}\Res_{z_1}Y_M(v,z_2)Y_M(u,z_1)w\\
&&\ \ \ \ \ \ \qquad\qquad\qquad\times\frac{(1+z_1)^{\wt u+l}}{z_1^{l+m+1}}\cdot\frac{(1+z_2)^{\wt v+l-1+\delta_i(r)+\frac{r}{T}}}{z_2^{2l+\delta_i(r)+\delta_i(T-r)+1}}\\
&&+\sum_{m=0}^l(-1)^m{m+l\choose l}\sum_{j,k\geq 0}{\wt u+l\choose j}{-l-m-1\choose k}\Res_{z_2}Y_M(u_{j+k}v,z_2)w\\
&&\ \ \ \ \ \ \qquad\qquad\qquad\times\frac{(1+z_2)^{\wt u+\wt v+2l-j-1+\delta_i(r)+\frac{r}{T}}}{z_2^{3l+\delta_i(r)+\delta_i(T-r)+2+m+k}}\\
&&=\sum_{m=0}^l(-1)^m{m+l\choose l}\Res_{z_2}\Res_{z_1}Y_M(v,z_2)Y_M(u,z_1)w\\
&&\ \ \ \ \ \ \qquad\qquad\qquad\times\frac{(1+z_1)^{\wt u+l}}{z_1^{l+m+1}}\cdot\frac{(1+z_2)^{\wt v+l-1+\delta_i(r)+\frac{r}{T}}}{z_2^{2l+\delta_i(r)+\delta_i(T-r)+1}}\\
&&+\sum_{m=0}^l(-1)^m{m+l\choose l}\sum_{k,j\geq 0}{\wt u+l\choose j}{-l-m-1\choose k}\Res_{z_2}Y_M(u_{k+j}v,z_2)w\\
&&\ \ \ \ \ \ \qquad\qquad\qquad\times\frac{(1+z_2)^{\wt (u_{k+j}v)+l-1+\delta_i(r)+\frac{r}{T}+l+k+1}}{z_2^{2l+\delta_i(r)+\delta_i(T-r)+1+l+m+k+1}}
\end{eqnarray*}
By the definition of $\O_{g,n}(M)$ and Lemma \ref{l2.2}$(1)$, we know that the two resulting terms are in $\O_{g,n}(M)$.

Using Lemma \ref{l2.2}$(2)$, similarly, we also have $(v\circ_{g,n}w )*_{g,n}u\in\O_{g,n}(M).$
so the proof of the Lemma is complete. \qed
\end{pf}

Now we can establish the main results of this section as follows:
\begin{thm}\label{t3.5} Let $M$ be an admissible $V$-module and  $n\in\frac{1}{T}\Z_+.$ Then
\begin{itemize}
\item[(1)] The vector space $\AA_{g,n}(M)$ is an $A_{g,n}(V)$-bimodule with the left and right actions of $A_{g,n}(V)$ on $\AA_{g,n}(M)$
induced from $(\ref{2.2m})$ and $(\ref{2.3m})$ respectively.

\item[(2)] The identity map on $M$ induces an $A_{g,n}(V)$-bimodule epimorphism from $\AA_{g,n}(M)$ to $\AA_{g,n-\frac{1}{T}}(M)$
if $n\geq 1.$

\item[(3)] The map
$$\phi:  w\mapsto e^{L(1)}(-1)^{L(0)}w$$
induces a linear isomorphism from $\AA_{g,n}(M)$ to $\AA_{g^{-1},n}(M)$ such that
$$\phi(u*_{g,n}w)=\phi(w)*_{g^{-1},n}\phi(u), \phi(w*_{g,n}u)=\phi(u)*_{g^{-1},n}\phi(w)$$
for $u\in V$ and $w\in M.$

\item[(4)] If $V$ is $g$-rational, then both $\O_{g,n-\frac{s}{T}}(M)/\O_{g,n-\frac{s-1}{T}}(M)$ and $\AA_{g,n-\frac{s}{T}}(M)$ are $A_{g,n}(V)$-bimodules
for $s=1,...,nT$ and
$$\AA_{g,n}(M)=\AA_{g,0}(M)\bigoplus \bigoplus_{s=1}^{nT}\O_{g,n-\frac{s}{T}}(M)/\O_{g,n-\frac{s-1}{T}}(M).$$

\end{itemize}
\end{thm}

\pf $(1)$ By Lemmas \ref{l2.f1} and \ref{l2.3} it suffices to show that the following identities hold in $\AA_{g,n}(M)$ for $u,v\in V$ and $w\in M:$
\begin{eqnarray*}
& &(u*_{g,n}w)*_{g,n}v=u*_{g,n}(w*_{g,n}v), \1*_{g,n}w=w*_{g,n}\1=w\\
& & (u*_{g,n}v)*_{g,n}w=u*_{g,n}(v*_{g,n}w), w*_{g,n}(u*_{g,n}v)=(w*_{g,n}u)*_{g,n}v.
\end{eqnarray*}
As the proofs are similar, we only prove the third identity in detail. From the definition, we may assume $u, v\in V^0$. Then we have
\begin{eqnarray*}
& &(u*_{g,n}v)*_{g,n}w=\sum_{m_1,m_2=0}^l\sum_{i\geq 0}(-1)^{m_1+m_2}{m_1+l\choose
l}{m_2+l\choose l}{\wt u+l\choose i}\\
& &\cdot \Res_{z_2}Y_M(u_{i-l-m_1-1}v,z_2)w\frac{(1+z_2)^{\wt u+\wt v+2l+m_1-i}}{z_2^{l+m_2+1}}\\
& &=\sum_{m_1,m_2=0}^l(-1)^{m_1+m_2}{m_1+l\choose
l}{m_2+l\choose l}\\
& &\cdot \Res_{z_2}\Res_{z_0}Y_M(Y(u,z_0)v,z_2)w\frac{(1+z_2+z_0)^{\wt\,u+l}(1+z_2)^{\wt v+m_1+l}}{z_0^{l+m_1+1}z_2^{l+m_2+1}}\\
& &=\sum_{m_1,m_2=0}^l(-1)^{m_1+m_2}{m_1+l\choose
l}{m_2+l\choose l}\\
& &\cdot \Res_{z_1}\Res_{z_2}Y_M(u,z_1)Y_M(v,z_2)w\frac{(1+z_1)^{\wt\,u+l}(1+z_2)^{\wt v+m_1+l}}{(z_1-z_2)^{l+m_1+1}z_2^{l+m_2+1}}\\
& &-\sum_{m_1,m_2=0}^l(-1)^{m_1+m_2}{m_1+l\choose
l}{m_2+l\choose l}\\
& &\cdot \Res_{z_2}\Res_{z_1}Y_M(v,z_2)Y_M(u,z_1)w\frac{(1+z_1)^{\wt\,u+l}(1+z_2)^{\wt v+m_1+l}}{(-z_2+z_1)^{l+m_1+1}z_2^{l+m_2+1}}\\
& &=\sum_{m_1,m_2=0}^l\sum_{i\geq 0}(-1)^{m_1+m_2+i}{m_1+l\choose
l}{m_2+l\choose l}{-l-m_1-1\choose i}\\
& &\cdot \Res_{z_1}\Res_{z_2}Y_M(u,z_{1})Y_M(v,z_2)w\frac{(1+z_1)^{\wt u+l}(1+z_2)^{\wt v+m_1+l}}{z_1^{l+m_1+1+i}z_2^{l+m_2+1-i}}\\
& &-\sum_{m_1,m_2=0}^l\sum_{i\geq 0}(-1)^{m_2+l+1+i}{m_1+l\choose
l}{m_2+l\choose l}{-l-m_1-1\choose i}\\
& &\cdot \Res_{z_2}\Res_{z_1}Y_M(v,z_2)Y_M(u,z_1)w\frac{(1+z_1)^{\wt u+l}(1+z_2)^{\wt v+m_1+l}}{z_1^{-i}z_2^{2l+m_1+m_2+2+i}}.
\end{eqnarray*}
It follows from
Lemma \ref{l2.2}$(1)$ that $$\Res_{z_2}\Res_{z_1}Y_M(v,z_2)Y_M(u,z_1)w\frac{(1+z_1)^{\wt u+l}(1+z_2)^{\wt v+l+m_1}}{z_1^{-i}z_2^{2l+m_1+m_2+2+i}}\in\O_{g,n}(M).$$ and if $i>l-m_1$
$$\Res_{z_1}\Res_{z_2}Y_M(u,z_{1})Y_M(v,z_2)w\frac{(1+z_1)^{\wt u+l}(1+z_2)^{\wt v+m_1+l}}{z_1^{l+m_1+1+i}z_2^{l+m_2+1-i}}$$ is also in $\O_{g,n}(M).$ Thus in $\AA_{g,n}(M)$, we have
\begin{eqnarray*}
& &(u*_{g,n}v)*_{g,n}w=u*_{g,n}(v*_{g,n}w)+
\sum_{m_1,m_2=0}^l(-1)^{m_1+m_2}{m_1+l\choose
l}{m_2+l\choose l}\\
& &\cdot \Res_{z_1}\Res_{z_2}Y_M(u,z_{1})Y_M(v,z_2)w\frac{(1+z_1)^{\wt u+l}(1+z_2)^{\wt v+l}}{z_1^{l+m_1+1}z_2^{l+m_2+1}}\\
& &\cdot (\sum_{i=0}^{l-m_1}\sum_{j\geq 0}{-l-m_1-1\choose i}{m_1\choose j}\frac{(-1)^iz_2^{i+j}}{z_1^i}-1)
\end{eqnarray*}
By proposition 5.3 in \cite{DLM3}, we know that $$\sum_{m_1=0}^l(-1)^{m_1}{m_1+l\choose
l}(\sum_{i=0}^{l-m_1}\sum_{j\geq 0}{-l-m_1-1\choose i}{m_1\choose j}\frac{(-1)^iz_2^{i+j}}{z_1^{i+m_1}}-\frac{1}{z_1^{m_1}})=0.$$
So we obtain $$(u*_{g,n}v)*_{g,n}w=u*_{g,n}(v*_{g,n}w).$$

$(2)$ By definition and Lemma \ref{l2.2}$(1)$, it is easy to see that $\O_{g,n}(M)\subset \O_{g,n-\frac{1}{T}}(M).$ So it suffices to show
that $u*_{g,n}w\equiv u*_{g,n-\frac{1}{T}}w, w*_{g,n}u\equiv w*_{g,n-\frac{1}{T}}u$ modulo  $\O_{g,n-\frac{1}{T}}(M)$ for $u\in V$ and $w\in M.$
We only need to prove them for $u\in V^0.$

As $n=l+\frac{i}{T}, n-\frac{1}{T}=l+\frac{i-1}{T}$, if $i\geq 1$, for $n$ and $n-\frac{1}{T}$, the nonnegative integer $l$ stays the same. It is clear that $u*_{g,n}w=u*_{g,n-\frac{1}{T}}w, w*_{g,n}u=w*_{g,n-\frac{1}{T}}u$. If $i=0$, then $n=l, n-\frac{1}{T}=l-1+\frac{T-1}{T}$. By definition
(\ref{2.2m}) we have

\begin{eqnarray*}
& &(u*_{g,n}w)=\sum_{m=0}^l(-1)^{m}{m+l\choose
l}\Res_{z}Y_M(u,z)w\frac{(1+z)^{\wt u+l}}{z^{l+m+1}}\\
& &=\sum_{m=0}^l(-1)^{m}{m+l\choose
l}\Res_{z}Y_M(u,z)w(\frac{(1+z)^{\wt u+l-1}}{z^{l+m}}+\frac{(1+z)^{\wt u+l-1}}{z^{l+m+1}})\\
& &=\sum_{m=0}^{l-1}(-1)^{m}{m+l\choose
l}\Res_{z}Y_M(u,z)w\frac{(1+z)^{\wt u+l-1}}{z^{l+m}}+\\
& &\sum_{m=0}^{l-2}(-1)^{m}{m+l\choose
l}\Res_{z}Y_M(u,z)w\frac{(1+z)^{\wt u+l-1}}{z^{l+m+1}}\quad (\mod\O_{g,n-\frac{1}{T}}(M))\\\end{eqnarray*}
\begin{eqnarray*}
& &=\Res_{z}Y_M(u,z)w\frac{(1+z)^{\wt u+l-1}}{z^{l}}+\sum_{m=1}^{l-1}(-1)^{m}{m+l\choose
l}\Res_{z}Y_M(u,z)w\frac{(1+z)^{\wt u+l-1}}{z^{l+m}}\\
& &+\sum_{m=1}^{l-1}(-1)^{m-1}{m+l-1\choose
l}\Res_{z}Y_M(u,z)w\frac{(1+z)^{\wt u+l-1}}{z^{l+m}}\\
& &=\Res_{z}Y_M(u,z)w\frac{(1+z)^{\wt u+l-1}}{z^{l}}+\sum_{m=1}^{l-1}\Res_{z}Y_M(u,z)w\frac{(1+z)^{\wt u+l-1}}{z^{l+m}}\\
& &\cdot ((-1)^{m}{m+l\choose l}+(-1)^{m+1}{m+l-1\choose l})\\
& &=\Res_{z}Y_M(u,z)w\frac{(1+z)^{\wt u+l-1}}{z^{l}}+\sum_{m=1}^{l-1}(-1)^{m}{m+l-1\choose
l-1}\Res_{z}Y_M(u,z)w\frac{(1+z)^{\wt u+l-1}}{z^{l+m}}\\
& &=\sum_{m=0}^{l-1}(-1)^{m}{m+l-1\choose
l-1}\Res_{z}Y_M(u,z)w\frac{(1+z)^{\wt u+l-1}}{z^{l+m}}\\
& &=u*_{g,n-\frac{1}{T}}w\quad (\mod\O_{g,n-\frac{1}{T}}(M)).
\end{eqnarray*}
Using this result and Lemma \ref{l2.2}, this complete the proof of $(2)$.

$(3)$ We fist prove that $\phi(\O_{g,n}(M))\subset \O_{g^{-1},n}(M).$  We know the following conjugation formulas from \cite{FHL}:
$$z^{L(0)}Y(u,z_{0})z^{-L(0)}=Y(z^{L(0)}u,zz_{0}),$$
$$e^{zL(1)}Y(u,z_{0})e^{-zL(1)}=Y\left(e^{z(1-zz_{0})L(1)}(1-zz_{0})^{-2L(0)}u,{z_0\over 1-zz_{0}}\right)$$
on $M$ for $u\in V.$ Then for homogeneous $u\in V^r$ and $w\in M,$
\begin{eqnarray*}
& &\phi(u\circ_{g,n} w)=e^{L(1)}(-1)^{L(0)}\Res_{z}\frac{(1+z)^{{\wt}u+l-1+\delta_i(r)+\frac{r}{T}}}{z^{2l+\delta_i(r)+\delta_i(T-r)+1}}Y_M(u,z)w
\end{eqnarray*}
\begin{eqnarray*}
& &=\Res_{z}\frac{(1+z)^{{\wt}u+l-1+\delta_i(r)+\frac{r}{T}}}{z^{2l+\delta_i(r)+\delta_i(T-r)+1}}e^{L(1)}
Y_M((-1)^{L(0)}u,-z)(-1)^{L(0)}w\\
& &=\Res_{z}\frac{(1+z)^{{\wt}u+l-1+\delta_i(r)+\frac{r}{T}}}{z^{2l+\delta_i(r)+\delta_i(T-r)+1}}
Y_M(e^{(1+z)L(1)}(1+z)^{-2L(0)}(-1)^{L(0)}u,\frac{-z}{1+z})e^{L(1)}(-1)^{L(0)}w.
\end{eqnarray*}
Setting variable $z=\displaystyle{-{z_{0}\over 1+z_{0}}}$ and
using the residue formula for the change of variable \cite{Z}
$$\Res_zg(z)=\Res_{z_0}(g(f(z_0)){d\over{dz_0}}f(z_0))$$
(for $g(z)=\sum_{n\geq N}v_nz^{n+\alpha}$ and $f(z)=\sum_{n>0}c_nz^n$
with $n\in\Z$, $v_n\in V,$ $\alpha, c_n\in\C$)

we have
\begin{eqnarray*}
& &\ \ \ \ \phi(u\circ_{g,n}w)\\
& &=(-1)^{\delta_i(r)+\delta_i(T-r)}\Res_{z_{0}}\frac{(1+z_0)^{-{\wt}u+l-1+\delta_i(T-r)+\frac{T-r}{T}}}{z_0^{2l+\delta_i(r)+\delta_i(T-r)+1}}\\
& &\times Y_M\left(e^{(1+z_{0})^{-1}L(1)}(1+z_{0})^{2L(0)}(-1)^{L(0)}u,z_{0}\right)
 e^{L(1)}(-1)^{L(0)}w\\
& &=(-1)^{{\wt}u+\delta_i(r)+\delta_i(T-r)}\Res_{z_0}\frac{(1+z_0)^{{\wt}u+l-1+\delta_i(T-r)+\frac{T-r}{T}}}{z_0^{2l+\delta_i(r)+\delta_i(T-r)+1}}\\
& &\times Y_M(e^{(1+z_0)^{-1}L(1)}u,z_0)e^{L(1)}(-1)^{L(0)}w\nonumber\\
& &=(-1)^{{\wt}u+\delta_i(r)+\delta_i(T-r)}\sum_{j=0}^{\infty}{1\over j!}\Res_{z_0}
\frac{(1+z_0)^{{\wt}u+l-1+\delta_i(T-r)+\frac{T-r}{T}}-j}{z_0^{2l+\delta_i(r)+\delta_i(T-r)+1}}\\
& &\times Y_M(L(1)^{j}u,z)e^{L(1)}(-1)^{L(0)}w\\
& &=(-1)^{{\wt}u+\delta_i(r)+\delta_i(T-r)}\sum_{j=0}^{\infty}{1\over j!}\Res_{z_0}
\frac{(1+z_0)^{{\wt}(L(1)^ju)+l-1+\delta_i(T-r)+\frac{T-r}{T}}}{z_0^{2l+\delta_i(r)+\delta_i(T-r)+1}}\\
& &\times Y_M(L(1)^{j}u,z_0)e^{L(1)}(-1)^{L(0)}w
\end{eqnarray*}
which lies in $\O_{g^{-1},n}(M)$ by definition.

Next we prove $\phi(u*_{g,n}w)=\phi(w)*_{g^{-1},n}\phi(u)$ and $\phi(w*_{g,n}u)=\phi(u)*_{g^{-1},n}\phi(w)$.

For $u\in V^r, r\neq 0$, from the definition (\ref{2.2m}) it is clear that $$\phi(u*_{g,n}w)=\phi(w)*_{g^{-1},n}\phi(u)=0.$$
Now assume $u\in V^0$. Then we have
\begin{eqnarray*}
&&\phi(u*_{g,n}w)=\phi\left(\sum_{m=0}^{l}(-1)^m{m+l\choose l}\Res_zY_M(u,z)w\frac{(1+z)^{\wt\,u+l}}{z^{l+m+1}}\right)\\
& &=\sum_{m=0}^{l}(-1)^m{m+l\choose l}\Res_z\frac{(1+z)^{\wt\,u+l}}{z^{l+m+1}}e^{L(1)}Y_M((-1)^{L(0)}u,-z)(-1)^{L(0)}w\\
& &=\sum_{m=0}^{l}(-1)^m{m+l\choose l}\Res_z\frac{(1+z)^{\wt\,u+l}}{z^{l+m+1}}\cdot
Y(e^{(1+z)L(1)}(1+z)^{-2L(0)}(-1)^{L(0)}u,\frac{-z}{1+z})e^{L(1)}(-1)^{L(0)}w\\
& &=\sum_{m=0}^{l}(-1)^{\wt u+l}{m+l\choose l}\Res_{z_0}\frac{(1+z_0)^{\wt\,u+m-1}}{z_0^{l+m+1}}Y_M(e^{(1+z_0)^{-1}L(1)}u,z_0)e^{L(1)}(-1)^{L(0)}w\\
& &=\sum_{j=0}^{\infty}\frac{1}{
j!}\sum_{m=0}^{l}(-1)^{\wt u+l}{m+l\choose l}\Res_{z_0}\frac{(1+z_0)^{\wt\,u+m-1-j}}{z_0^{l+m+1}}
Y(L(1)^{j}u,z_0)e^{L(1)}(-1)^{L(0)}w\\
& &=\phi(w)*_{g^{-1},n}\sum_{j=0}^{\infty}\frac{1}{
j!}(L(1)^j(-1)^{L(0)}u)\\
& &=\phi(w)*_{g^{-1},n}\phi(u).
\end{eqnarray*}
where the forth step uses change of variable $z=\displaystyle{-{z_{0}\over 1+z_{0}}}$.

The proof for the other identity is similar, so we omit it.

$(4)$ From Lemma \ref{l2.f1} and \ref{l2.3}, it is obvious that $\O_{g,n-\frac{s}{T}}(M)/\O_{g,n-\frac{s-1}{T}}(M)$ and $\AA_{g,n-\frac{s}{T}}(M)$ are the $A_{g,n}(V)$-bimodules for $s=1,2,\cdots, nT.$ By $(2)$ we have $$\psi:  \AA_{g,n}(M)\longrightarrow \AA_{g,n-\frac{1}{T}}(M)$$ a algebra epimorphism. Then
$$ \AA_{g,n-\frac{1}{T}}(M)\cong \AA_{g,n}(M)/\ker\psi.$$
Since $\O_{g,n}(M)\subset \O_{g,n-\frac{1}{T}}(M)$, we know that $$\ker\psi =\O_{g,n-\frac{1}{T}}(M)/\O_{g,n}(M).$$
If $V$ is $g$-rational, then
\begin{eqnarray*}\AA_{g,n}(M)&&=\AA_{g,n-\frac{1}{T}}(M)\bigoplus\O_{g,n-\frac{1}{T}}(M)/\O_{g,n}(M)\\
&&=\AA_{g,n-\frac{2}{T}}(M)\bigoplus\O_{g,n-\frac{2}{T}}(M)/\O_{g,n-\frac{1}{T}}(M)\bigoplus\O_{g,n-\frac{1}{T}}(M)/\O_{g,n}(M)\\
&&=\cdots\cdots\cdots\\
&&=\AA_{g,0}(M)\bigoplus\bigoplus_{s=1}^{nT}\O_{g,n-\frac{s}{T}}(M)/\O_{g,n-\frac{s-1}{T}}(M)
\end{eqnarray*}\qed

\section{$\AA_{g,n}(M)$ and intertwining operators}
\hspace{0.5cm}
Now we discuss the connections between $\AA_{g,n}(M)$ and the intertwining operators and fusion rules.Throughout this section, let $(M_{0},Y_{M_0})$ be an admissible $V$-module and $(M_{i},Y_{M_i})$ be admissible $g$-twisted V-modules for $i=1,2$.
We first recall the intertwining operators and tensor product of modules from \cite{DLM7}.

\begin{definition}An intertwining operator of type \hspace{-0.2
cm}\singlespace $\left(\!\hspace{-3 pt}\begin{array}{c} M_2\\
M_0 M_1\end{array}\hspace{-3 pt}\!\right)$\doublespace is a linear
map
\begin{eqnarray*}
M_{0}& &\longrightarrow ({\rm Hom}(M_{1},M_{2}))\{z\}\\
w & &\mapsto I(w,z)=\sum_{n\in\mathbb{C}}w_{n}z^{-n-1}
\end{eqnarray*}
such that for $w\in M_0, w^1\in M_1$, fixed $c\in\mathbb{C}$ and $n\in\Q$ sufficiently large
$$w_{n+c}w^1=0,$$ the following (generalized) Jacobi identity holds on $M_2$:
for $u\in V^r$ and $w^j\in M_j, j=0,1,$
\begin{eqnarray*}
&& z_{0}^{-1}\delta
(\frac{z_{1}-z_{2}}{z_{0}})Y_{M_2}(u,z_{1})I(w^0,z_{2})w^1-
z_{0}^{-1}\delta (\frac{z_{2}-z_{1}}{-z_{0}})I(w^0,z_{2})Y_{M_1}(u,z_{1})w^1\\
&=& z_{2}^{-1}(\frac{z_{1}-z_{0}}{z_{2}})^{-\frac{r}{T}}\delta
(\frac{z_{1}-z_{0}}{z_{2}})I(Y_{M_0}(u,z_{0})w^0,z_{2})w^1
\end{eqnarray*}
and $$I(L(-1))w^0,z)=\frac{d}{dz}I(w^0,z).$$
\end{definition}
Remark: Assume $M_0$ is an irreducible $V$-module, and $M_1$ $M_2$ are irreducible $g$-twisted $V$-modules. We use $h_i\in \mathbb{C}, ~i=0,1,2$ to denote the conformal weights(ref. Remark \ref{rem1}) of $M_i,~ i=0,1,2.$ Due to Proposition 1.5.1 of \cite{FZ}, the intertwining operator of type \hspace{-0.2
cm}\singlespace $\left(\!\hspace{-3 pt}\begin{array}{c} M_2\\
M_0 M_1\end{array}\hspace{-3 pt}\!\right)$ can be written as
$$I(w,z)=\sum_{n\in\frac{1}{T}\Z_+}w_n z^{-n-1}z^{-h_0-h_1+h_2},
$$
so that for homogeneous $w\in M_0,$ $m,n\in \frac{1}{T}\Z_+,$
\begin{eqnarray*}
w(n) M_{1}(m)\subseteq M_{2}(\deg w+m-n-1),
\end{eqnarray*}
 where $\deg w=k$ means that $w\in M_0(k).$

The intertwining operators of type \hspace{-0.2
cm}\singlespace $\left(\!\hspace{-3 pt}\begin{array}{c} M_2\\
M_0 M_1\end{array}\hspace{-3 pt}\!\right)$\doublespace clearly form a vector space, which we denote by ${\cal I}\left (\hspace{-3 pt}\begin{array}{c} M_2\\
M_0\,M_1\end{array}\hspace{-3 pt}\right)$. Set
$$N_{M_{0}\,M_{1}}^{M_{2}}=\dim\ \hspace{-0.2
cm}\singlespace {\cal I}\left (\hspace{-3 pt}\begin{array}{c} M_2\\
M_0\,M_1\end{array}\hspace{-3 pt}\right)$$
which are called the fusion rules. For short, we set $N_{01}^2=N_{M_{0}\,M_{1}}^{M_{2}}.$
We also set ${\cal I}_{01}^2={\cal I}\left (\hspace{-3 pt}\begin{array}{c} M_2\\
M_0\,M_1\end{array}\hspace{-3 pt}\right).$

Now we can define the tensor product of an admissible $V$-modules and an admissible $g$-twisted $V$-module according to \cite{DLM7}.
\begin{definition} A tensor product for the admissible module $M_0$ and the admissible $g$-twisted module $M_1$ is an admissible $g$-twisted module $M_0\boxtimes M_1$ together with an intertwining operator
$F\in {\cal I}\left(\hspace{-3 pt}\begin{array}{c} M_0\boxtimes M_1\\
M_0\,M_1\end{array}\hspace{-3 pt}\right)$\doublespace
such that the following universal mapping property holds:
for any a $g$-twisted $V$-module $W$ and any
intertwining operator $I\in {\cal I}\left(\hspace{-3 pt}\begin{array}{c} W\\
M_0\,M_1\end{array}\hspace{-3 pt}\right),$\doublespace there exists a unique
$V$-homomorphism $\psi$ from $M_0\boxtimes M_1$ to $W$ such that $I=\psi\circ F.$
\end{definition}

From now on we assume that $V$ is $g$-rational. Let $M^1,...,M^p$ be the irreducible admissible $g$-twisted $V$-modules up to equivalence.
The following result is well known:
\begin{thm}\label{tensor} If $V$ is rational and $g$-rational then the tensor product $M_0\boxtimes M^j$ exists for any $j.$
In fact,
$$M_{0}\boxtimes M_j=\bigoplus_{k=1}^p N_{0j}^kM^{k}$$
such that  $I_{0j}^{k}(w,z)M^j\subset M^{k}\{z\}$ for $w\in M_0.$
\end{thm}

We now turn our attention to the connection between intertwining operators and $\AA_{g,n}(M)$ as before.
From now on, we denote the untwisted admissible module $M$ by $M_0$.
Let
$I\in {\cal I}_{0j}^k.$ Then for $w\in M_0$,
$w(\deg w-1-t+s)M^j(s)\subset M^k(t)$ for all $s,t.$ For short we set
$$o_{t,s}^I(w)=w(\deg w-1-t+s)$$
 for homogeneous $w\in M_0$ and extend it linearly to entire $M_0.$ According to the Lemma 1.5.2 of \cite{FZ} and Theorem 3.2 of \cite{DLM4} and Lemma 4.3 of \cite{DR}, we have the following lemma.
\begin{lem}\label{ll}   Let $M_0, M_1,M_2$ be as before, which are all irreducible and $I$ is an intertwining operator of type \hspace{-0.2
cm}\singlespace $\left(\!\hspace{-3 pt}\begin{array}{c} M_2\\
M_0 M_1\end{array}\hspace{-3 pt}\!\right).$The map
\begin{equation*}
\begin{split}
\pi(I):\AA_{g,n}(M_0)\otimes_{A_{g,n}(V)}M_1(s) &\to M_2(t)\\
(w^0,w^1)&\mapsto o^I_{t,s}(w^0)w^1
\end{split}
\end{equation*}
for $w^0\in M_0, w^1\in M_1(s)$ and $s,t\in\frac{1}{T}\Z_+,$ such that $s,t \leq n$ induces
an $A_{g,n}(V)$-module homomorphism.
\end{lem}

\pf As in \cite{FZ} we need to verify that $o^I_{t,s}(u*_{g,n}w^0)=o(u)o^I_{t,s}(w^0),$ $o^I_{t,s}(w^0*_{g,n}u)=o^I_{t,s}(w^0)o(u)$ and
$o^I_{t,s}(w)=0$ on $M_1(s)$ with $s,t\leq n$ for $u\in V^0, w^0\in M_0$ and $w\in \O_{g,n}(M_0).$

By the definition, we have
\begin{eqnarray*}& &o^I_{t,s}(u*_{g,n}w^0)=o^I_{t,s}(\sum_{m=0}^l(-1)^m{m+l \choose l}\Res_zY_{M_0}(u,z)w^0\frac{(1+z)^{\wt u+l}}{z^{l+m+1}})\\
& &=o^I_{t,s}(\sum_{m=0}^l\sum_{i\geq 0}(-1)^m{m+l \choose l}{\wt u+l \choose i}u_{i-l-m-1}w^0)\\
& &=\sum_{m=0}^l\sum_{i\geq 0}(-1)^m{m+l \choose l}{\wt u+l \choose i}(u_{i-l-m-1}w^0)_{\wt u+\deg w^0-i+l+m-1-t+s}\\\end{eqnarray*}
\begin{eqnarray*}
& &=\sum_{m=0}^l(-1)^m{m+l \choose l}\Res_{z_2}\Res_{z_0}I_{M_1}(Y_{M_0}(u,z_0)w^0,z_2)\times\\
& & \frac{(z_2+z_0)^{\wt u+l}z_2^{h_0+h_1-h_2+\deg w^0+m-1-t+s}}{z_0^{l+m+1}}\\
& &=\sum_{m=0}^l(-1)^m{m+l \choose l}\Res_{z_1}\Res_{z_2}Y_{M_2}(u,z_1)I(w^0,z_2)\times\\
& & \frac{z_1^{\wt u+l}z_2^{h_0+h_1-h_2+\deg w^0+m-1-t+s}}{(z_1-z_2)^{l+m+1}}\\
& &-\sum_{m=0}^l(-1)^m{m+l \choose l}\Res_{z_2}\Res_{z_1}I(w^0,z_2)Y_{M_1}(u,z_1)\times\\
& & \frac{z_1^{\wt u+l}z_2^{h_0+h_1-h_2+\deg w^0+m-1-t+s}}{(-z_2+z_1)^{l+m+1}}\\
& &=\sum_{m=0}^l\sum_{j\geq 0}(-1)^{m+j}{m+l \choose l}{-l-m-1 \choose j}u_{\wt u-m-1-j}w^0_{\deg w^0+m-1-t+s+j}\\
& &-\sum_{m=0}^l\sum_{j\geq 0}(-1)^{l+m+1+j}{m+l \choose l}{-l-m-1 \choose j}w^0_{\deg w^0-l-2-t+s-j}u_{\wt u+l+j}\\
& &=\sum_{k=0}^\infty\sum_{m=0}^k(-1)^{k}{m+l \choose l}{-l-m-1 \choose k-m}u_{\wt u-1-k}w^0_{\deg w^0-1-t+s+k}\\
& &-\sum_{j=0}^\infty\sum_{m=0}^j(-1)^{l+1+j}{m+l \choose l}{-l-m-1 \choose j}w^0_{\deg w^0-l-2-t+s-j}u_{\wt u+l+j}
\end{eqnarray*}
Note that $u_{\wt u+l+j}=w^0_{\wt w^0-1-t+s+k}=0$ on $M_1(s)$ for $j\geq 0$ and $k>t.$ Also we have
$$\sum_{m=0}^k{m+l \choose l}{-l-m-1 \choose k-m}=\sum_{m=0}^k{m+l \choose m}{-l-m-1 \choose k-m}=0, \text{for}\ 1\leq k\leq l.$$
Then $$o^I_{t,s}(u*_{g,n}w^0)=u_{\wt u-1}w^0_{\wt w^0-1-t+s}=o(u)o^I_{t,s}(w^0).$$
Similarly, one can get
$$o^I_{t,s}(w^0*_{g,n}u)=w^0_{\deg w^0-1-t+s}u_{\wt u-1}=o^I_{t,s}(w^0)o(u).$$
Next we prove the third identity. For $w\in\O_{g,n}(M_0),$ setting $w=a\circ_{g,n}b$ for $a\in V^r, b\in M_0$, we have
\begin{eqnarray*}& &o^I_{t,s}(a\circ_{g,n}b)=o^I_{t,s}(\Res_zY_{M_0}(a,z)b\frac{(1+z)^{\wt a+l-1+\delta_i(r)+\frac{r}{T}}}{z^{2l+\delta _i(r)+\delta_i(T-r)+1}})\\\end{eqnarray*}
\begin{eqnarray*}
& &=o^I_{t,s}(\sum_{i\geq 0}{\wt a+l-1+\delta_i(r)+\frac{r}{T} \choose i}a_{i-2l-1-\delta _i(r)-\delta_i(T-r)}b)\\
& &=\sum_{i\geq 0}{\wt a+l-1+\delta_i(r)+\frac{r}{T} \choose i}(a_{i-2l-1-\delta _i(r)-\delta_i(T-r)}b)_{\wt u+\deg b-i+2l+\delta _i(r)+\delta_i(T-r)-1-t+s}\\
& &=\Res_{z_2}\Res_{z_0}I(Y_{M_0}(a,z_0)b,z_2)\frac{(z_2+z_0)^{\wt a+l-1+\delta_i(r)+\frac{r}{T}}z_2^{h_0+h_1-h_2+\deg b+l+\delta_i(T-r)-t+s-\frac{r}{T}}}{z_0^{2l+\delta _i(r)+\delta_i(T-r)+1}}\\
& &=\Res_{z_1}\Res_{z_2}Y_{M_2}(a,z_1)I(b,z_2)\frac{z_1^{\wt a+l-1+\delta_i(r)+\frac{r}{T}}z_2^{h_0+h_1-h_2+\deg b+l+\delta_i(T-r)-t+s-\frac{r}{T}}}{(z_1-z_2)^{2l+\delta _i(r)+\delta_i(T-r)+1}}\\
& &-\Res_{z_2}\Res_{z_1}I(b,z_2)Y_{M_1}(a,z_1)\frac{z_1^{\wt a+l-1+\delta_i(r)+\frac{r}{T}}z_2^{h_0+h_1-h_2+\deg b+l+\delta_i(T-r)-t+s-\frac{r}{T}}}{(-z_2+z_1)^{2l+\delta _i(r)+\delta_i(T-r)+1}}\\
& &=\sum_{j\geq 0}(-1)^j{-2l-\delta_i(r)-\delta_i(T-r)-1 \choose j}\Res_{z_1}\Res_{z_2}Y_{M_2}(a,z_1)I(b,z_2)\times\\
& &\frac{z_1^{\wt a+l-1+\delta_i(r)+\frac{r}{T}}z_2^{h_0+h_1-h_2+\deg b+l+\delta_i(T-r)-t+s-\frac{r}{T}}}{z_1^{2l+\delta _i(r)+\delta_i(T-r)+1+j}z_2^{-j}}\\
& &-\sum_{j\geq 0}(-1)^{\delta_i(r)+\delta_i(T-r)+1}{-2l-\delta_i(r)-\delta_i(T-r)-1 \choose j}\Res_{z_2}\Res_{z_1}I(b,z_2)Y_{M_1}(a,z_1)\times\\
& &\frac{z_1^{\wt a+l-1+\delta_i(r)+\frac{r}{T}}z_2^{h_0+h_1-h_2+\deg b+l+\delta_i(T-r)-t+s-\frac{r}{T}}}{z_1^{-j}z_2^{2l+\delta _i(r)+\delta_i(T-r)+1+j}}\\
& &=\sum_{j\geq 0}(-1)^j{-2l-\delta_i(r)-\delta_i(T-r)-1 \choose j}
a_{\wt a-l-j-2-\delta_i(T-r)+\frac{r}{T}}b_{\deg b+l+\delta _i(T-r)-t+s-\frac{r}{T}+j}\\
& &+\sum_{j\geq 0}(-1)^{\delta_i(r)+\delta_i(T-r)}{-2l-\delta_i(r)-\delta_i(T-r)-1 \choose j}
b_{\deg b-l-j-\delta_i(r)-1-\frac{r}{T}-t+s}a_{\wt a+l-1+\delta_i(r)+\frac{r}{T}+j}\end{eqnarray*}
Note that $$b_{\deg b+l+\delta _i(T-r)-t+s-\frac{r}{T}+j}=b_{\deg b+s-1+l+\delta _i(T-r)-t+\frac{T-r}{T}+j}=0$$ and $$a_{\wt a+l-1+\delta_i(r)+\frac{r}{T}+j}=0$$ on $M_1(s)$ for $j\geq 0$.
So $o^I_{t,s}(w)=0$.
This completes the proof of Lemma.
\qed

We now obtain a linear map $\pi$ from $ {\cal I}\left(\hspace{-3 pt}\begin{array}{c} M_2\\
M_0\,M_1\end{array}\hspace{-3 pt}\right)$ to $\Hom_{A_{g,n}(V)}(\AA_{g,n}(M_0)\otimes_{A_{g,n}(V)}M_1(s), M_2(t)).$ The following lemma is from Proposition 2.10. of \cite{L2}.

\begin{lem}\label{12} If $M_2$ is irreducible, $\pi$ is injective.\qed\end{lem}

By Lemmas \ref{ll} and \ref{12}, a similar argument as that in \cite{L2} for the proof of Theorem 1.5.2 of \cite{FZ} and Theorem 2.11 in \cite{L2} could be applied to our twisted case, which gives us  the main result of this section:
\begin{thm}\label{tn1} Suppose $M_0$ is an irreducible $V$-module, $M_1$ and $M_2$ are irreducible $g$-twisted $V$-modules, Then $\pi$ is an isomorphism from $ {\cal I}\left(\hspace{-3 pt}\begin{array}{c} M_2\\
M_0\,M_1\end{array}\hspace{-3 pt}\right)$ to $\Hom_{A_{g,n}(V)}(\AA_{g,n}(M_0)\otimes_{A_{g,n}(V)}M_1(s), M_2(t))$ for $s,t\leq n$.
\end{thm}
\begin{rem}From the above, in general, we obtain that there is an isomorphism from $ {\cal I}\left(\hspace{-3 pt}\begin{array}{c} M^k\\
M_0\,M^j\end{array}\hspace{-3 pt}\right)$ to $\Hom_{A_{g,n}(V)}(\AA_{g,n}(M_0)\otimes_{A_{g,n}(V)}M^j(s), M^k(t))$ for $s,t\leq n$ where $M^j, M^k$ are irreducible $g$-twisted $V$-module.\end{rem}

\end{document}